\documentclass[12pt]{article}

\usepackage{amsmath}
\allowdisplaybreaks[4]

\usepackage[all]{xy}
\usepackage{amssymb}
\usepackage{amsthm}
\usepackage{hyperref}
\hypersetup{colorlinks=true,linkcolor=blue,citecolor=red}
\usepackage{amsmath}
\usepackage{amscd}
\usepackage{verbatim}
\usepackage{eurosym}
\usepackage{float}
\usepackage{color}
\usepackage{dcolumn}
\usepackage[mathscr]{eucal}
\usepackage[all]{xy}
\usepackage{hyperref}
\usepackage{mathrsfs}
\usepackage{amsmath}
\usepackage{amssymb}
\usepackage{amsfonts,ifpdf}
\usepackage{graphicx}
\usepackage{times}
\usepackage{float}
\usepackage{epstopdf}
\usepackage{cite}
\usepackage{youngtab}
\usepackage{ytableau}
\ytableausetup
{mathmode, boxsize=0.9em}

\newlength{\boxedparwidth}
{\begin{center} \begin{tabular}{|@{\hspace{.315in}}c@{\hspace{.15in}}|}
			\hline \\ \begin{minipage}[t]{\boxedparwidth}
				\setlength{\parindent}{.25in}}%
			{\end{minipage} \\ \\ \hline \end{tabular} \end{center}}

\allowdisplaybreaks

\newtheorem{thm}{Theorem}[section]
\newtheorem{prop}[thm]{Proposition}

\newtheorem{lem}[thm]{Lemma}
\newtheorem{cor}[thm]{Corollary}

\numberwithin{equation}{section}

\begin{document}
		
\begin{center}		
{\Large \bf Extensions of $q$-supercongruences arising from a quadratic summation of Gasper and Rahman}
\end{center}
\vskip 3mm

\begin{center}
		
Qun Mu and Ying Zhong\\[8pt]
 School of Mathematics\\
Hunan University\\
Changsha 410082, P.R. China\\[12pt]

Emails:  qunmu@hnu.edu.cn,  YingZhong@hnu.edu.cn
\\[15pt]		
		
\end{center}
	
\vskip 3mm
	
\begin{abstract}
		
\vskip 6pt
In this paper, inspired by recent work of Guo and Zhao, and of Guo and Zudilin, we establish two parametric extensions of their $q$-supercongruences using Gasper and Rahman's quadratic summation, the Chinese remainder theorem, and the creative microscoping method. These extensions yield two one-parameter generalizations of the $q$-supercongruences considered by Guo and Zhao and by Guo and Zudilin.

\vskip 6pt
		
\\ [7pt]
	\end{abstract}
	
\section{Introduction}
In 1914, Ramanujan \cite{Ramanujan-1914} established the remarkable identity
	\begin{equation}\label{R1}
			\sum_{k=0}^{\infty}(6k+1) \frac{\left(\frac{1}{2}\right)_{k}^{3}}{k!^{3} 4^{k}}=\frac{4}{\pi},
	\end{equation}
	where $(a)_{n}=a(a+1)\dots (a+n-1)$ is the Pochhammer symbol. Motivated by Ramanujan's series for $\frac{1}{\pi}$, Van Hamme \cite{VanHamme-1997} conjectured $13$ $p$-adic analogues, now known as Ramanujan-type supercongruences. In particular, the conjecture associated with \eqref{R1} is that for primes $p > 3$,
	\begin{equation}\label{(J.2)}
		\sum_{k=0}^{\frac{p-1}{2}}
		(6k+1)\frac{\left(\frac12\right)_k^3}
		{k!^3\,4^k}
		\equiv
		(-1)^{\frac{p-1}{2}}p
		\pmod{p^4},
	\end{equation}
which was marked as (J.2) in his list. Early progress on \eqref{(J.2)} was obtained by Zudilin \cite{Zudilin-2009}, who verified the supercongruence modulo $p^2$ using the WZ method. The first complete proof of \eqref{(J.2)} was given by Long \cite{Long-2011}, using hypergeometric identities and $p$-adic analysis. Later, Long and Ramakrishna \cite{Long-Ramakrishna-2016} established a related supercongruence modulo $p^6$ by means of $p$-adic Gamma function and hypergeometric series formulas.
Further similar Ramanujan-type supercongruences have been established by
various authors; see, for example, \cite{He-2017, Liu-2017}.


Building on these results, the $q$-analogues of Ramanujan's and Ramanujan-type $\pi$-series have attracted considerable attention. Guo and Liu \cite{Guo-Liu-2018} derived a $q$-analogue of the Ramanujan's $\pi$-series \eqref{R1}. For further references on this topic, see \cite{Chen-Chu-2021, Chu-2018,Hou-K-Sun-2019,Sun-2019}. Meanwhile, Ramanujan-type supercongruences have been generalized to the $q$-settings by many authors. Guo \cite{Guo-2018} first established a weaker $q$-analogue of the supercongruence \eqref{(J.2)}, and subsequently proved the complete $q$-analogue in \cite{Guo-2020-1} using the $q$-microscope method.
This method was introduced by Guo and Zudilin \cite{Guo-Zudilin-2019} and has subsequently led to numerous extensions of $q$-supercongruences; see, for example, \cite{Guo-2020, Guo-2022, Guo-2024, Gu-Guo-2022, Tang-2023, Liu-Wang-2022, Wang-2024, Wang-Xu-2025, Wei-2022, Wei-2024}.

Recently, Guo and Zudilin \cite[Theorem 4.5]{Guo-Zudilin-2019} showed that
\begin{align}\label{cong:GZ}
&\sum_{k=0}^{m}[6k+1]\frac{(aq;q^{2})_{k} (q/a;q^{2})_{k} (q;q^{2})_{k} (b;q^{4})_{k} }{(aq^{4};q^{4})_{k} (q^{4}/a;q^{4})_{k} (q^{4};q^{4})_{k} (q^{3}/b;q^{2})_{k} b^{k}} q^{k^{2}+2k }
\notag\\
&\qquad\equiv \frac{(q^{r}b;q^{4})_{(n-r)/4}}{ (q^{4+r}/b;q^{4})_{(n-r)/4}}b^{-(n-r)/4}(-q)^{(1-r)/2}[n] \pmod{[n](1-aq^{n})(a-q^{n})},
\end{align}
where $m=n-1$ or $\frac{n-1}2$, and $n \equiv r  \pmod  {4}$ be an odd positive integer with $r=\pm 1$. Here and throughout the paper, we adopt the standard $q$-notations.
For complex numbers $a$, $q$, and any nonnegative integer $n$, the $q$-shifted factorials are defined as
\[(a;q)_\infty = \prod_{j=0}^{\infty}(1 - aq^j)\quad \text{and} \quad (a;q)_n = \prod_{j=0}^{n-1}(1 - aq^j). \]
For simplicity, we shall also adopt the compact notation
\[(a_1, a_2, \dots, a_\ell; q)_n = (a_1;q)_n(a_2;q)_n\cdots(a_\ell;q)_n\]
for $n\in\mathbb{N}_0\cup\{\infty\}$.
In addition, the $q$-integer is given by $[n] = (1 - q^n)/(1 - q) $.
Let $\Phi_n(q)$ denote the $n$-th cyclotomic polynomial, which can be written as
		\begin{equation*}
		\begin{aligned}
			\Phi_{n}(q)=\prod_{\substack{j=1 \\ \gcd(j,n)=1}}^{n-1} (q-e^{2\pi ij/n})\in \mathbb{Z}[q] .
		\end{aligned}
	\end{equation*}
Letting $a \to 1$, $b \to q^{-1}$ and $r=-1$ in \eqref{cong:GZ}, we have
	\begin{equation}\label{GZu1}	
			\sum_{k=0}^{m}[6k+1]\frac{(q;q^2)_{k}^{3} (q^{-1};q^{4})_{k}}{ (q^{4};q^{2})_{k} (q^{4};q^{4})_{k}^{3}}q^{k^2+3k} \equiv -[n]q^{(n+5)/4} \frac{(q^{-2};q^{4})_{(n+1)/4}}{(q^{4};q^{4})_{(n+1)/4}}  \pmod{[n]\Phi_{n}(q)^{2}}.
	\end{equation}
More recently, Guo and Zhao
\cite[Theorem 1.3]{Guo-Zhao-2025} gave a related $q$-supercongruences modulo $[n]\Phi _{n}(q) ^{3}$, that is, for any positive integer $n$ with $n \equiv 3  \pmod  {4}$,
	\begin{equation}\label{GZh}
			\sum_{k=0}^{\frac{n-1}{2} }\left [ 6k+1 \right ] \frac{  \left ( q;q^{2} \right )_{k}^{3}     \left ( q^{-1};q^{4} \right )_{k}  }{  \left ( q^{4};q^{4} \right )_{k}^{3}     \left ( q^{4};q^{2} \right )_{k}  }q^{k^{2}+3k} \equiv -[n]q^{\frac{n+5}{4}  }\frac{\left(q^{-2} ; q^{4}\right)_{(n+1) / 4}}{\left(q^{4} ; q^{4}\right)_{(n+1) / 4}} \left [ 1-[n]^{2}  \sum_{j=1}^{(n+1) / 4} \frac{q^{4 j}}{[4 j]^{2}} \right ]	.		
	\end{equation}

In this paper, we aim to investigate parametric extensions of the $q$-supercongruences \eqref{GZu1} and \eqref{GZh} under different moduli. More precisely, we first establish two $q$-supercongruences modulo $\Phi _{n}(q) ^{3}$ and $\Phi _{n}(q) ^{4}$, respectively, and then, building on these results, derive one-parameter generalizations of the $q$-supercongruences \eqref{GZu1} and \eqref{GZh}.
	

We first establish an extension of \eqref{GZu1} as follows.

\begin{thm}\label{thm:mod-phi3}
		Let $n> 1$ be an integer with $n \equiv 3\ (\bmod \ 4)$, and $s$ be an integer with $0\le s \le \frac{n-1}{6}$. Then, modulo $\Phi _{n}(q) ^{3}$,
		\begin{equation*}
			\begin{aligned}
				&\sum_{k=s}^{\frac{n-1}{2}+2s } \left [ 6k+1 \right ] \frac{ \left ( q;q^{2} \right )_{k}  \left ( q;q^{2} \right )_{k-2s}  \left ( q;q^{2} \right )_{k+2s}  \left ( q^{-1};q^{4} \right )_{k}  }{ \left ( q^{4};q^{4} \right )_{k+s}  \left ( q^{4};q^{4} \right )_{k}  \left ( q^{4};q^{4} \right )_{k-s}  \left ( q^{4};q^{2} \right )_{k}  }q^{k^{2}+3k} \\
				& \quad \equiv  -[n]q^{\frac{n+5}{4}  }\frac{\left(q^{-2} ; q^{4}\right)_{(n+1) / 4}}{\left(q^{4} ; q^{4}\right)_{(n+1) / 4}} \frac{ \left ( q;q^{2} \right ) _{s} \left ( q;q^{2} \right ) _{-s} \left ( q^{-1};q^{4} \right ) _{s} }{    \left ( q^{5} ; q^{4}\right ) _{ s } }q^{s^{2}+3s}.
			\end{aligned}
		\end{equation*}
	\end{thm}
The following result provides a one-parameter extension of \eqref{GZh} modulo $\Phi _{n}(q) ^{4}$.

\begin{thm}\label{thm:analogue}
Let $n> 1$ be an integer with $n \equiv 3\ (\bmod \ 4)$, and $s$ be an integer with $0\le s \le  \frac{n-1}{6}$. Then, modulo $\Phi _{n}(q) ^{4}$,
\begin{equation*}
\begin{aligned}
&\sum_{k=s}^{\frac{n-1}{2}+2s } \left [ 6k+1 \right ] \frac{ \left ( q;q^{2} \right )_{k}  \left ( q;q^{2} \right )_{k-2s}  \left ( q;q^{2} \right )_{k+2s}  \left ( q^{-1};q^{4} \right )_{k}  }{ \left ( q^{4};q^{4} \right )_{k+s}  \left ( q^{4};q^{4} \right )_{k}  \left ( q^{4};q^{4} \right )_{k-s}  \left ( q^{4};q^{2} \right )_{k}  }q^{k^{2}+3k}
\\
& \quad \equiv -[n]q^{\frac{n+5}{4}  }\frac{\left(q^{-2} ; q^{4}\right)_{(n+1) / 4}}{\left(q^{4} ; q^{4}\right)_{(n+1) / 4}} \frac{ \left ( q;q^{2} \right ) _{s} \left ( q;q^{2} \right ) _{-s} \left ( q^{-1};q^{4} \right ) _{s} }{    \left ( q^{5} ; q^{4}\right ) _{ s } } \left [ 1-[n]^{2}  \sum_{j=1}^{(n+1) / 4} \frac{q^{4 j}}{[4 j]^{2}} \right ] q^{s^{2}+3s} .
\end{aligned}
\end{equation*}
\end{thm}

Next, from Theorem \ref{thm:mod-phi3}, we derive a generalization of \eqref{GZu1}.
\begin{thm}\label{thm.1.4}
Let $n> 1$ be an integer with $n \equiv 3\ (\bmod \ 4)$, and $s$ be an integer with $0\le s \le \frac{n-1}{6}$. Then, modulo $[n]\Phi _{n}(q) ^{2}$,
\begin{align}\label{them 4}
&\sum_{k=s}^{\frac{n-1}{2}+2s } \left [ 6k+1 \right ] \frac{ \left ( q;q^{2} \right )_{k}  \left ( q;q^{2} \right )_{k-2s}  \left ( q;q^{2} \right )_{k+2s}  \left ( q^{-1};q^{4} \right )_{k}  }{ \left ( q^{4};q^{4} \right )_{k+s}  \left ( q^{4};q^{4} \right )_{k}  \left ( q^{4};q^{4} \right )_{k-s}  \left ( q^{4};q^{2} \right )_{k}  }q^{k^{2}+3k}
\notag\\
& \quad \equiv -[n]q^{\frac{n+5}{4}  }\frac{\left(q^{-2} ; q^{4}\right)_{(n+1) / 4}}{\left(q^{4} ; q^{4}\right)_{(n+1) / 4}} \frac{ \left ( q;q^{2} \right ) _{s} \left ( q;q^{2} \right ) _{-s} \left ( q^{-1};q^{4} \right ) _{s} }{    \left ( q^{5} ; q^{4}\right ) _{ s } }q^{s^{2}+3s} .
\end{align}
\end{thm}

Setting $s = 0$ in \eqref{them 4} yields \eqref{GZu1}. Moreover, letting $n = p$ be an odd prime and taking the limit $q \to 1$ in \eqref{them 4}, we obtain the following supercongruence.
\begin{cor}	
Let $p>1$ be an odd prime satisfying $p\equiv 3\pmod{4}$,
and let $s$ be an integer with $0\le s \le \frac{p-1}{6}$. Then we have
\begin{equation*}
\sum_{k=s}^{\frac{(p-1)}{2}+2s }(6k+1)\frac{\left ( \frac{1}{2} \right ) _{k} \left ( \frac{1}{2} \right ) _{k-2s} \left ( \frac{1}{2} \right )  _{k+2s} \left ( -\frac{1}{4} \right )  _{k}}{4^{k}(k+s)!k!(k-s)! (k+1)!}
				\equiv -p\frac{(-1)^{s}\left ( -\frac{1}{2} \right ) _{\frac{(p+1)}{4} } \left ( -\frac{1}{4} \right ) _{s}}{\left ( 1 \right ) _{\frac{(p+1)}{4} } \left ( \frac{5}{4} \right ) _{s }} \pmod{p^3}.		\end{equation*}
\end{cor}
	
Finally, with the aid of Theorem \ref{thm:analogue}, we obtain a generalization of \eqref{GZh}.
	\begin{thm}\label{thm.1.6}
		Let $n> 1$ be an integer with $n \equiv 3\ (\bmod \ 4)$, and $s$ be an integer with $0\le s \le  \frac{n-1}{6}$. Then, modulo $[n]\Phi _{n}(q) ^{3}$,
		\begin{equation}\label{them 6}
			\begin{aligned}
				&\sum_{k=s}^{\frac{n-1}{2}+2s } \left [ 6k+1 \right ] \frac{ \left ( q;q^{2} \right )_{k}  \left ( q;q^{2} \right )_{k-2s}  \left ( q;q^{2} \right )_{k+2s}  \left ( q^{-1};q^{4} \right )_{k}  }{ \left ( q^{4};q^{4} \right )_{k+s}  \left ( q^{4};q^{4} \right )_{k}  \left ( q^{4};q^{4} \right )_{k-s}  \left ( q^{4};q^{2} \right )_{k}  }q^{k^{2}+3k} \\
				&\quad \equiv  -[n]q^{\frac{n+5}{4}  }\frac{\left(q^{-2} ; q^{4}\right)_{(n+1) / 4}}{\left(q^{4} ; q^{4}\right)_{(n+1) / 4}} \frac{ \left ( q;q^{2} \right ) _{s} \left ( q;q^{2} \right ) _{-s} \left ( q^{-1};q^{4} \right ) _{s} }{    \left ( q^{5} ; q^{4}\right ) _{ s } } \left [ 1-[n]^{2}  \sum_{j=1}^{(n+1) / 4} \frac{q^{4 j}}{[4 j]^{2}} \right ] q^{s^{2}+3s} .
			\end{aligned}
		\end{equation}
	\end{thm}
	
	Putting $s = 0$ in \eqref{them 6} gives \eqref{GZh}. Furthermore, for $n = p$ an odd prime, taking the limit $q \to 1$ in \eqref{them 6} yields the following supercongruence.
	\begin{cor}	
Let $p>1$ be an odd prime satisfying $p\equiv 3\pmod{4}$,
and let $s$ be an integer satisfying $0\le s \le \frac{p-1}6$. Then we have
\begin{align*}
&\sum_{k=s}^{\frac{(p-1)}{2}+2s }(6k+1)\frac{\left ( \frac{1}{2} \right ) _{k} \left ( \frac{1}{2} \right ) _{k-2s} \left ( \frac{1}{2} \right )  _{k+2s} \left ( -\frac{1}{4} \right )  _{k}}{4^{k}(k+s)!k! (k-s)!(k+1)!}
\\
&\qquad\qquad\equiv -p\frac{(-1)^{s}\left ( -\frac{1}{2} \right ) _{\frac{(p+1)}{4} } \left ( -\frac{1}{4} \right ) _{s}}{\left ( 1 \right ) _{\frac{(p+1)}{4} } \left ( \frac{5}{4} \right ) _{s }}
				\left(1-\frac{p^2}{16}\sum_{j=1}^{(p+1)/4}\frac1{j^2}\right) \pmod{p^4}.		\end{align*}
	\end{cor}

The rest of the paper is organized as follows. In Section \ref{Sect.2} we give some necessary lemmas on $q$-congruences. In Section \ref{Sect.3}, we prove Theorems \ref{thm:mod-phi3} and \ref{thm:analogue} by the creative microscoping method. Finally, Section \ref{Sect.4} is devoted to the proofs of Theorems \ref{thm.1.4} and \ref{thm.1.6}.
	
\section{Some Lemmas}\label{Sect.2}	
In order to prove our main results, we establish four lemmas concerning $q$-congruences. To this end, we need the following formula
\begin{align}\label{equ:1.2}
\sum_{k=0}^{\infty}  \frac{1-a q^{3 k}}{1-a} \frac{(a, b, q / b ; q)_{k}\left(f ; q^{2}\right)_{k} q^{\left(k^{2}+k\right) / 2}}{\left(q^{2}, a q^{2} / b, a b q ; q^{2}\right)_{k}(a q / f ; q)_{k}}\left(\frac{a}{f}\right)^{k}
=\frac{\left(a q, a q^{2}, a q^{2} / b f, a b q / f ; q^{2}\right)_{\infty}}{\left(a q / f, a q^{2} / f, a q^{2} / b, a b q ; q^{2}\right)_{\infty}},
\end{align}
which follows by putting $d=q^{-2n}$ and taking $n\rightarrow \infty $ in the quadratic summation of Gasper and Rahman \cite[Eq. (3.8.12)]{Gasper-Rahman-2004}
\begin{align*}
\sum_{k=0}^{\infty} & \frac{1-a q^{3 k}}{1-a} \frac{(a, b, q / b ; q)_{k}\left(d, f, a^{2} q / d f ; q^{2}\right)_{k} q^{k}}{(a q / d, a q / f, d f / a ; q)_{k}\left(q^{2}, a q^{2} / b, a b q ; q^{2}\right)_{k}} \\
+ & \frac{(a q, f / a, b, q / b ; q)_{\infty}\left(d, a^{2} q / d f, f q^{2} / d, d f^{2} q / a^{2} ; q^{2}\right)_{\infty}}{(a / f, f q / a, a q / d, d f / a ; q)_{\infty}\left(a q^{2} / b, a b q, f q / a b, b f / a ; q^{2}\right)_{\infty}}
\\
\times& \sum_{k=0}^{\infty} \frac{\left(f, b f / a, f q / a b ; q^{2}\right)_{k} q^{2 k}}{\left(q^{2}, f q^{2} / d, d f^{2} q / a^{2} ; q^{2}\right)_{k}}
\\
=\ & \frac{(a q, f / a ; q)_{\infty}\left(a q^{2} / b d, a b q / d, b d f / a, d f q / a b ; q^{2}\right)_{\infty}}{(a q / d, d f / a ; q)_{\infty}\left(a q^{2} / b, a b q, b f / a, f q / a b ; q^{2}\right)_{\infty}}.
\end{align*}
Define
\begin{equation}\label{def:T}
T_{n,s}(a,b;q):=\sum_{k=s}^{\frac{n-1}{2}+2s }\left [ 6k+1 \right ] \frac{ \left ( aq;q^{2} \right )_{k}  \left ( q/a;q^{2} \right )_{k-2s}  \left ( q;q^{2} \right )_{k+2s}  \left ( q^{-1}/b;q^{4} \right )_{k}  }{ \left ( aq^{4};q^{4} \right )_{k+s}  \left ( q^{4}/a;q^{4} \right )_{k}  \left ( q^{4};q^{4} \right )_{k-s}  \left ( bq^{4};q^{2} \right )_{k}  }b^{k}q^{k^{2}+3k}.
\end{equation}

\begin{lem}\label{lem 1}
Let $n> 1$ be an odd integer, and $s$ be an integer with $0\le s \le \frac{n-1}{6}$. Then
\begin{equation}\label{equ:phi}
T_{n,s}(a,b;q) \equiv 0 \pmod{\Phi _{n}\left ( q \right )}.
\end{equation}
	\end{lem}
	
\begin{proof}
Replacing $k$ by $k+s$, a direct computation gives that $T_{n,s}(a,b;q)$ is equal to
\begin{align}\label{eq:2.2}
&\sum_{k=0}^{\frac{n-1}{2}+s } \left [ 6k+6s+1 \right ] \frac{ \left ( aq;q^{2} \right )_{k+s}  \left ( q/a;q^{2} \right )_{k-s}  \left ( q;q^{2} \right )_{k+3s}  \left ( q^{-1}/b;q^{4} \right )_{k+s}  }{ \left ( aq^{4};q^{4} \right )_{k+2s}  \left ( q^{4}/a;q^{4} \right )_{k+s}  \left ( q^{4};q^{4} \right )_{k}  \left ( bq^{4};q^{2} \right )_{k+s}  }q^{(k+s)^{2}+3(k+s)}b^{k+s}
\notag\\
&\quad  = \left [ 6s+1 \right ] \frac{ \left ( aq;q^{2} \right ) _{s} \left ( q/a;q^{2} \right ) _{-s} \left ( q;q^{2} \right ) _{3s} \left ( q^{-1}/b;q^{4} \right ) _{s}}{\left ( aq^{4};q^{4} \right ) _{2s} \left ( q^{4}/a;q^{4} \right ) _{s} \left ( bq^{4};q^{2} \right ) _{s} }b^{s}q^{s^{2}+3s}
\notag\\
&\qquad  \times \sum_{k=0}^{\frac{n-1}{2}+s } \frac{1-q^{6k+6s+1}}{1-q^{6s+1}}  \frac{ \left ( aq^{2s+1} , q^{1-2s}/a , q^{6s+1};q^{2} \right )_{k}  \left ( q^{4s-1}/b ;q^{4} \right )_{k}    }{ \left ( aq^{8s+4} , q^{4s+4}/a , q^{4};q^{4} \right )_{k}  \left ( bq^{2s+4};q^{2} \right )_{k} }b^{k}q^{k^{2}+2ks+3k}.
\end{align}

If $s = \frac{n-1}{6}$, then the right-hand side of \eqref{eq:2.2} contains
$[6s+1]=[n]$, whereas the summand equals $1$ for $k=0$ and
\[\frac{(q^{6s+1};q^{2})_{k}}{1-q^{6s+1}}=(q^{n+2};q^{2})_{k-1},\quad k\geq 1.\]
On the other hand, $(q^4;q^4)_k$ is relatively prime to
$\Phi_n(q)$ since $0\le k\le \frac{n-1}{2}+s$. Hence the desired congruence follows.

Now suppose that $0\le s < \frac{n-1}{6}$. Applying \eqref{equ:1.2} with $q\mapsto q^{2}$, $a=q^{6s+1-n}$, $b= aq^{2s+1}$ and $f=q^{4s-1}/b$, and along with the fact that $( q^{6s+1-n};q^2)_k=0$ for $k>\frac{n-1}{2}+s $, we have
		\begin{align*}
			&\sum_{k=0}^{\frac{n-1}{2}+s} \frac{1-q^{6k+6s+1-n}}{1-q^{6s+1-n}}  \frac{ \left ( aq^{2s+1} , q^{1-2s}/a , q^{6s+1-n};q^{2} \right )_{k}  \left ( q^{4s-1}/b ;q^{4} \right )_{k}    }{ \left ( aq^{8s+4-n} , q^{4s+4-n}/a , q^{4};q^{4} \right )_{k}  \left ( bq^{2s+4-n};q^{2} \right )_{k} }b^{k}q^{k^{2}+2ks+3k-nk}\\
			&\quad =\frac{\left ( q^{6s+3-n} , q^{6s+5-n} , bq^{5-n}/a , abq^{4s+5-n} ; q^{4}\right ) _{\infty } }{\left ( bq^{2s+4-n} , bq^{2s+6-n} , q^{4s+4-n}/a , aq^{8s+4-n} ; q^{4}\right ) _{\infty }}=0.
		\end{align*}
The last equality follows from
\[\left ( q^{6s+3-n} , q^{6s+5-n};q^4\right ) _{\infty }=
\left ( q^{6s+3-n};q^2\right ) _{\infty }=0,\]
because $6s+3-n$ is a nonpositive even integer.
Using $q^n\equiv1\pmod{\Phi_n(q)}$, we may therefore remove each
shift by $-n$ in the parameters and replace $q^{-nk}$ by $1$.
It follows that modulo $\Phi_{n}(q)$,
\begin{align*}
&\sum_{k=0}^{\frac{n-1}{2}+s}
\frac{1-q^{6k+6s+1}}{1-q^{6s+1}}
\frac{(aq^{2s+1},q^{1-2s}/a,q^{6s+1};q^2)_k
      (q^{4s-1}/b;q^4)_k}
     {(aq^{8s+4},q^{4s+4}/a,q^4;q^4)_k
      (bq^{2s+4};q^2)_k}
\,b^kq^{k^2+2sk+3k}
\equiv0.
\end{align*}
Plugging this into \eqref{eq:2.2}, we complete the proof.
\end{proof}

\begin{lem}\label{lem 2}
Let $n> 1$ be a positive integer with $n \equiv 3\ (\bmod \ 4)$, and $s$ be a non-negative integer. Then, modulo $a-q^{n}$,
\begin{align}\label{equ:a-q^{n}}
T_{n,s}(a,b;q) \equiv \frac{ \left ( aq;q^{2} \right ) _{s} \left ( q/a;q^{2} \right ) _{-s} \left ( q^{-1}/b;q^{4} \right ) _{s} \left ( q^{3};q^{4} \right ) _{\frac{n+1}{4}} \left (  bq^{5}/a ; q^{4}\right ) _{\frac{n+1}{4}} }{  \left ( q^{4}/a;q^{4} \right ) _{\frac{n+1}{4}} \left ( bq^{4} ; q^{4}\right ) _{\frac{n+1}{4} } \left ( abq^{5} ; q^{4}\right ) _{ s }}b^{s}q^{s^{2}+3s}.
\end{align}
\end{lem}
	
\begin{proof}
Letting $a=q^{n}$ in \eqref{def:T}, we obtain that $T_{n,s}(q^{n},b;q)$ is equal to
		\begin{align}\label{eq:2.7}
			&\sum_{k=s}^{\frac{n-1}{2}+2s } \left [ 6k+1 \right ] \frac{ \left ( q^{1+n};q^{2} \right )_{k}  \left ( q^{1-n};q^{2} \right )_{k-2s}  \left ( q;q^{2} \right )_{k+2s}  \left ( q^{-1}/b;q^{4} \right )_{k}  }{ \left ( q^{4+n};q^{4} \right )_{k+s}  \left ( q^{4-n};q^{4} \right )_{k}  \left ( q^{4};q^{4} \right )_{k-s}  \left ( bq^{4};q^{2} \right )_{k}  }q^{k^{2}+3k}b^{k}\nonumber\\
			& \quad = \left [ 6s+1 \right ] \frac{ \left ( q^{1+n};q^{2} \right ) _{s} \left ( q^{1-n};q^{2} \right ) _{-s} \left ( q;q^{2} \right ) _{3s} \left ( q^{-1}/b;q^{4} \right ) _{s}}{\left ( q^{4+n};q^{4} \right ) _{2s} \left ( q^{4-n};q^{4} \right ) _{s} \left ( bq^{4};q^{2} \right ) _{s} }b^{s}q^{s^{2}+3s} \nonumber\\
			&\qquad  \times \sum_{k=0}^{\frac{n-1}{2}+s } \frac{1-q^{6k+6s+1}}{1-q^{6s+1}}  \frac{ \left ( q^{2s+1+n} , q^{1-2s-n} , q^{6s+1};q^{2} \right )_{k}  \left ( q^{4s-1}/b ;q^{4} \right )_{k}    }{ \left ( q^{8s+4+n} , q^{4s+4-n} , q^{4};q^{4} \right )_{k}  \left ( bq^{2s+4};q^{2} \right )_{k} }b^{k}q^{k^{2}+2ks+3k}.
		\end{align}		
Setting $q\mapsto q^{2}$, $a=q^{6s+1}$, $b= q^{2s+n+1}$ and $f = q^{4s-1}/b$ in \eqref{equ:1.2}, we have
		\begin{align*}
			&\sum_{k=0}^{\frac{n-1}{2}+s } \frac{1-q^{6k+6s+1}}{1-q^{6s+1}}  \frac{ \left ( q^{2s+1+n} , q^{1-2s-n} , q^{6s+1};q^{2} \right )_{k}  \left ( q^{4s-1}/b ;q^{4} \right )_{k}    }{ \left ( q^{8s+4+n} , q^{4s+4-n} , q^{4};q^{4} \right )_{k}  \left ( bq^{2s+4};q^{2} \right )_{k} }b^{k}q^{k^{2}+2ks+3k}\\
			&\quad=\frac{\left ( q^{6s+3} , q^{6s+5} , bq^{5-n} , bq^{4s+n+5} ; q^{4}\right ) _{\infty } }{\left ( bq^{2s+4} , bq^{2s+6} , q^{4s-n+4} , q^{8s+n+4} ; q^{4}\right ) _{\infty }},
		\end{align*}
where we use that $( q^{1-2s-n};q^2)_k=0$ for $k>\frac{n-1}{2}+s $.
Then we can express the right-hand side of \eqref{eq:2.7} as
		\begin{align*}
		\frac{ \left ( q^{1+n};q^{2} \right ) _{s} \left ( q^{1-n};q^{2} \right ) _{-s} \left ( q^{-1}/b;q^{4} \right ) _{s} \left ( q^{3};q^{4} \right ) _{\frac{n+1}{4}  } \left (  bq^{5-n} ; q^{4}\right ) _{\frac{n+1}{4}} }{  \left ( q^{4-n};q^{4} \right ) _{\frac{n+1}{4}} \left ( bq^{4} ; q^{4}\right ) _{\frac{n+1}{4} } \left ( bq^{n+5} ; q^{4}\right ) _{ s }}b^{s}q^{s^{2}+3s}.
		\end{align*}
		This is the $a=q^{n}$ case of the right-hand side of \eqref{equ:a-q^{n}}. Namely, the $q$-congruence \eqref{equ:a-q^{n}} holds.
	\end{proof}

\begin{lem}\label{lem 3}
Let $n> 1$ be a positive integer with $n \equiv 3\ (\bmod \ 4)$, and $s$ be an integer with $0\le s \le \frac{n-1}{2}$.  Then, modulo $1-aq^{n}$,
\begin{align}\label{equ:1-aq^{n}}
T_{n,s}(a,b;q) \equiv \frac{ \left ( aq;q^{2} \right ) _{s} \left ( q/a;q^{2} \right ) _{-s} \left ( q^{-1}/b;q^{4} \right ) _{s} \left ( q^{3};q^{4} \right ) _{\frac{n+1}{4} } \left (  abq^{5} ; q^{4}\right ) _{\frac{n+1}{4}} }{ \left (  abq^{5} ; q^{4}\right ) _{s}  \left ( aq^{4};q^{4} \right ) _{\frac{n+1}{4}} \left ( bq^{4};q^{4} \right ) _{\frac{n+1}{4}} }b^{s}q^{s^{2}+3s}.
\end{align}
\end{lem}
\begin{proof}
Setting $a=q^{-n}$ in \eqref{def:T}, we have that $T_{n,s}(q^{-n},b;q)$ is
\begin{align}\label{eq:2.9}
&\sum_{k=s}^{\frac{n-1}{2}+2s}\left [ 6k+1 \right ] \frac{ \left ( q^{1-n};q^{2} \right )_{k}  \left ( q^{1+n};q^{2} \right )_{k-2s}  \left ( q;q^{2} \right )_{k+2s}  \left ( q^{-1}/b;q^{4} \right )_{k}  }{ \left ( q^{4-n};q^{4} \right )_{k+s}  \left ( q^{4+n};q^{4} \right )_{k}  \left ( q^{4};q^{4} \right )_{k-s}  \left ( bq^{4};q^{2} \right )_{k}  }q^{k^{2}+3k}b^{k}
\nonumber\\
&\quad= \left [ 6s+1 \right ] \frac{ \left ( q^{1-n};q^{2} \right ) _{s} \left ( q^{1+n};q^{2} \right ) _{-s} \left ( q;q^{2} \right ) _{3s} \left ( q^{-1}/b;q^{4} \right ) _{s}}{\left ( q^{4-n};q^{4} \right ) _{2s} \left ( q^{4+n};q^{4} \right ) _{s} \left ( bq^{4};q^{2} \right ) _{s} }b^{s}q^{s^{2}+3s}
\nonumber\\
&\qquad\times \sum_{k=0}^{\frac{n-1}{2}+s } \frac{1-q^{6k+6s+1}}{1-q^{6s+1}}  \frac{ \left ( q^{2s+1-n} , q^{1-2s+n} , q^{6s+1};q^{2} \right )_{k}  \left ( q^{4s-1}/b ;q^{4} \right )_{k}    }{ \left ( q^{8s+4-n} , q^{4s+4+n} , q^{4};q^{4} \right )_{k}  \left ( bq^{2s+4};q^{2} \right )_{k} }b^{k}q^{k^{2}+2ks+3k}.
\end{align}
Because $0\le s \le \frac{n-1}{2}$, one can verify that $( q^{2s+1-n};q^2)_k=0$ if $k>\frac{n-1}{2}-s$.		
Then applying \eqref{equ:1.2} with $q\mapsto q^{2}$, $a=q^{6s+1}$, $b= q^{2s-n+1}$ and $f=q^{4s-1}/b$, we have
		\begin{align*}
			&\sum_{k=0}^{\frac{n-1}{2}+s } \frac{1-q^{6k+6s+1}}{1-q^{6s+1}}  \frac{ \left ( q^{2s+1-n} , q^{1-2s+n} , q^{6s+1};q^{2} \right )_{k}  \left ( q^{4s-1}/b ;q^{4} \right )_{k}    }{ \left ( q^{8s+4-n} , q^{4s+4+n} , q^{4};q^{4} \right )_{k}  \left ( bq^{2s+4};q^{2} \right )_{k} }b^{k}q^{k^{2}+2ks+3k}\\
			&\quad =\frac{\left ( q^{6s+3} , q^{6s+5} , bq^{5+n} , bq^{4s-n+5} ; q^{4}\right ) _{\infty } }{\left ( bq^{2s+4} , bq^{2s+6} , q^{4s+n+4} , q^{8s-n+4} ; q^{4}\right ) _{\infty }}  .
		\end{align*}
Hence, the right-hand side of \eqref{eq:2.9} becomes
		\begin{align*}
		\frac{ \left ( q^{1-n};q^{2} \right ) _{s} \left ( q^{1+n};q^{2} \right ) _{-s} \left ( q^{-1}/b;q^{4} \right ) _{s} \left ( q^{3};q^{4} \right ) _{\frac{n+1}{4}}  \left ( bq^{-n+5} ; q^{4}\right ) _{\frac{n+1}{4} } }{\left ( bq^{-n+5} ; q^{4}\right ) _{ s } \left ( q^{4-n};q^{4} \right ) _{\frac{n+1}{4}} \left ( bq^{4};q^{4} \right ) _{\frac{n+1}{4}} } b^{s}q^{s^{2}+3s}.
		\end{align*}
This implies our desired result.
\end{proof}

\begin{lem}\label{lem 4}
Let $n> 1$ be a positive integer with $n \equiv 3\ (\bmod \ 4)$, and $s$ be an integer with $0\le s \le  \frac{n+1}{4}$. Then, modulo $b-q^{n}$,
\begin{align}\label{equ:b-q^{n}}
T_{n,s}(a,b;q)\equiv \frac{ \left ( aq;q^{2} \right ) _{s} \left ( q/a;q^{2} \right ) _{-s} \left ( q^{-1}/b;q^{4} \right ) _{s} \left ( q^{3};q^{4} \right ) _{\frac{n+1}{4} } \left ( q^{5};q^{4} \right ) _{\frac{n+1}{4}}    }{  \left ( aq^{4};q^{4} \right ) _{\frac{n+1}{4}}  \left ( abq^{5};q^{4} \right ) _{s} \left ( q^{4}/a;q^{4} \right ) _{\frac{n+1}{4}} }b^{s}q^{s^{2}+3s}.
\end{align}
\end{lem}
\begin{proof}
In view of \eqref{def:T}, we have that $T_{n,s}(a,q^{n};q)$ is equal to
\begin{align}\label{eq:2.11}
&\sum_{k=s}^{\frac{n-1}{2}+2s } \left [ 6k+1 \right ] \frac{ \left ( aq;q^{2} \right )_{k}  \left ( q/a;q^{2} \right )_{k-2s}  \left ( q;q^{2} \right )_{k+2s}  \left ( q^{-1-n};q^{4} \right )_{k}  }{ \left ( aq^{4};q^{4} \right )_{k+s}  \left ( q^{4}/a;q^{4} \right )_{k}  \left ( q^{4};q^{4} \right )_{k-s}  \left ( q^{4+n};q^{2} \right )_{k}  }q^{k^{2}+3k+nk}
\nonumber\\
& \quad = \left [ 6s+1 \right ] \frac{ \left ( aq;q^{2} \right ) _{s} \left ( q/a;q^{2} \right ) _{-s} \left ( q;q^{2} \right ) _{3s} \left ( q^{-1-n};q^{4} \right ) _{s}}{\left ( aq^{4};q^{4} \right ) _{2s} \left ( q^{4}/a;q^{4} \right ) _{s} \left ( q^{4+n};q^{2} \right ) _{s} }q^{s^{2}+3s+ns}
\nonumber\\
&\qquad  \times \sum_{k=0}^{\frac{n-1}{2}+s } \frac{1-q^{6k+6s+1}}{1-q^{6s+1}}  \frac{ \left ( aq^{2s+1} , q^{1-2s}/a , q^{6s+1};q^{2} \right )_{k}  \left ( q^{4s-1-n} ;q^{4} \right )_{k}    }{ \left ( aq^{8s+4} , q^{4s+4}/a , q^{4};q^{4} \right )_{k}  \left ( q^{2s+4+n};q^{2} \right )_{k} }q^{k^{2}+2ks+3k+nk}.
\end{align}
Since $0\le s \le  \frac{n+1}{4}$, we see that $\left ( q^{4s-1-n} ;q^{4} \right )_{k}=0$ when $k>\frac{n+1}{4}-s$.
Then performing the parameter substitutions $q\mapsto q^{2}$, $a=q^{6s+1}$, $b=aq^{2s+1}$ and $f=q^{4s-1-n}$ in \eqref{equ:1.2} gives
		\begin{align*}
			&\sum_{k=0}^{\frac{n-1}{2}+s } \frac{1-q^{6k+6s+1}}{1-q^{6s+1}}  \frac{ \left ( aq^{2s+1} , q^{1-2s}/a , q^{6s+1};q^{2} \right )_{k}  \left ( q^{4s-1-n} ;q^{4} \right )_{k}    }{ \left ( aq^{8s+4} , q^{4s+4}/a , q^{4};q^{4} \right )_{k}  \left ( q^{2s+4+n};q^{2} \right )_{k} }q^{k^{2}+2ks+3k+nk}\\
			&\quad =\frac{\left ( q^{6s+3} , q^{6s+5} , q^{5+n}/a , aq^{4s+n+5} ; q^{4}\right ) _{\infty } }{\left ( q^{2s+4+n} , q^{2s+6+n} , q^{4s+4}/a , aq^{8s+4} ; q^{4}\right ) _{\infty }} .
		\end{align*}
Consequently, the right-hand side of \eqref{eq:2.11} can be rewritten as
\begin{align*}
\frac{ \left ( aq;q^{2} \right ) _{s} \left ( q/a;q^{2} \right ) _{-s} \left ( q^{-1-n};q^{4} \right ) _{s}  \left ( q^{3};q^{4} \right ) _{\frac{n+1}{4} } \left ( q^{5};q^{4} \right ) _{\frac{n+1}{4} }  }{  \left ( aq^{4};q^{4} \right ) _{\frac{n+1}{4}} \left ( aq^{5+n};q^{4} \right ) _{s} \left ( q^{4}/a;q^{4} \right ) _{\frac{n+1}{4}} }q^{s^{2}+3s+ns},
\end{align*}
which completes the proof.
\end{proof}

\section{Proofs of Theorems \ref{thm:mod-phi3} and \ref{thm:analogue}}\label{Sect.3}
In this section, we prove Theorems \ref{thm:mod-phi3} and \ref{thm:analogue} by  the creative microscoping method. To this end, we provide two parametric extensions.
Building on the lemmas established in Section \ref{Sect.2} and the Chinese remainder theorem, we first derive the following parametric extension, which is used to prove Theorem \ref{thm:mod-phi3}.
	
\begin{prop}\label{prop:mod-phia}
Let $n> 1$ be a positive integer with $n \equiv 3\ (\bmod \ 4)$, and $s$ be an integer with $0\le s \le \frac{n-1}{6}$. Then, modulo $\Phi _{n}(q)(a-q^{n})(1-aq^{n})$,
\begin{align}\label{edu:mod all 3.1}
T_{n,s}(a,b;q)
&\equiv \frac{ \left ( aq;q^{2} \right ) _{s} \left ( q/a;q^{2} \right ) _{-s} \left ( q^{-1}/b;q^{4} \right ) _{s} \left ( q^{3};q^{4} \right ) _{\frac{n+1}{4}}  }{    \left ( abq^{5} ; q^{4}\right ) _{ s } \left ( bq^{4} ; q^{4}\right ) _{\frac{n+1}{4} }}
\notag\\
&\quad \times \left [ \frac{1-aq^{n}}{1-a^{2}} \frac{  \left (  bq^{5}/a ; q^{4}\right ) _{\frac{n+1}{4}} }{  \left ( q^{4}/a;q^{4} \right ) _{\frac{n+1}{4}}  }
				-\frac{a(a-q^{n})}{1-a^{2}}\frac{  \left (  abq^{5} ; q^{4}\right ) _{\frac{n+1}{4}} }{   \left ( aq^{4};q^{4} \right ) _{\frac{n+1}{4}} } \right ]b^{s}q^{s^{2}+3s}.
\end{align}

\end{prop}
	
\begin{proof}
Note that $ 1-aq^{n}$ and $ a-q^{n}$ are coprime polynomials. Applying the Chinese remainder theorem for coprime polynomials, together with Lemmas \ref{lem 2} and \ref{lem 3}, we can deduce the remainder of $T_{n,s}(a,b;q)$ modulo $(1-aq^{n})(a-q^{n}) $.  Thus, using
\begin{align*}
\frac{1-aq^{n}}{1-a^{2}} \equiv 1 \pmod{(a-q^{n})}
\quad \text{and} \quad
\frac{-a(a-q^{n})}{1-a^{2}} \equiv 1 \pmod{(1-aq^{n})} ,
\end{align*}
we conclude that \eqref{edu:mod all 3.1} holds modulo $(a-q^{n})(1-aq^{n})$.

Moreover, it is obvious that $\Phi_{n}(q)$ and $ (1-aq^{n})(a-q^{n})$ are coprime polynomials. Since the factor $\left ( q^{3};q^{4} \right ) _{\frac{n+1}{4}}$ on the right-hand side of \eqref{edu:mod all 3.1} contains $1-q^{n}$, the right-hand side of \eqref{edu:mod all 3.1} vanishes modulo $\Phi_{n}(q)$. By Lemma \ref{lem 1}, we can deduce that \eqref{edu:mod all 3.1} holds modulo $\Phi_{n}(q)(1-aq^{n})(a-q^{n})$.
This completes the proof.
\end{proof}
	
We next give a parametric congruence, which enables us to verify Theorem \ref{thm:analogue}.
	\begin{prop}\label{prop:mod all 3.3}
		Let $n> 1$ be a positive integer with $n \equiv 3\ (\bmod \ 4)$, and $s$ be an integer with $0\le s \le  \frac{n-1}{6}$. Then, modulo $\Phi _{n}(q)^{2}(a-q^{n})(1-aq^{n}) $,
\begin{align}\label{edu:mod all 3.3}
T_{n,s}(a,1;q)
&\equiv
\frac{(aq;q^2)_s(q/a;q^2)_{-s}(q^{-1};q^4)_s}
     {(aq^5;q^4)_s}q^{s^2+3s}
\notag\\
&\quad\times\left\{
\frac{(q^3,q^{5-n};q^4)_{\frac{n+1}{4}}}
     {(q^{4-n},q^4;q^4)_{\frac{n+1}{4}}}
+\frac{(q^3,q^5;q^4)_{\frac{n+1}{4}}}
      {(aq^4,q^4/a;q^4)_{\frac{n+1}{4}}}
-\frac{(q^3,q^5;q^4)_{\frac{n+1}{4}}}
      {(q^{4-n},q^{4+n};q^4)_{\frac{n+1}{4}}}
\right\}.
\end{align}
\end{prop}

\begin{proof}
We first prove the following two-parameter congruence, modulo
$\Phi_n(q)(a-q^n)(1-aq^n)(b-q^n)$:
\begin{align}\label{3.4}
T_{n,s}(a,b;q)
&\equiv
\frac{(aq;q^2)_s(q/a;q^2)_{-s}(q^{-1}/b;q^4)_s}
     {(abq^5;q^4)_s}b^s q^{s^2+3s}
\notag\\
&\quad\times\left\{
\frac{(q^3,bq^{5-n};q^4)_{\frac{n+1}{4}}}
     {(q^{4-n},bq^4;q^4)_{\frac{n+1}{4}}}
+\frac{(q^3,q^5;q^4)_{\frac{n+1}{4}}}
      {(aq^4,q^4/a;q^4)_{\frac{n+1}{4}}}
-\frac{(q^3,q^5;q^4)_{\frac{n+1}{4}}}
      {(q^{4-n},q^{4+n};q^4)_{\frac{n+1}{4}}}
\right\}.
\end{align}
For $a=q^n$ or $a=q^{-n}$, the second and third fractions in braces
cancel exactly. The remaining expression agrees with the corresponding
evaluation in Lemma \ref{lem 2} or Lemma \ref{lem 3}, respectively.
Thus \eqref{3.4} holds modulo $(a-q^n)(1-aq^n)$.
For $b=q^n$, the first and third fractions cancel exactly, and the
remaining expression is the one in Lemma \ref{lem 4}.
Hence \eqref{3.4} also holds modulo $b-q^n$.

Every fraction in braces contains $(q^3;q^4)_{\frac{n+1}{4}}$
in its numerator, and hence contains $1-q^n$.
Its denominator is relatively prime to $\Phi_n(q)$ for generic $a,b$.
Indeed, none of $4j-n$ and $4j+n$, with
$1\le j\le(n+1)/4<n$, is divisible by $n$.
By Lemma \ref{lem 1}, both sides of \eqref{3.4} therefore vanish
modulo $\Phi_n(q)$.
The four polynomial factors
$\Phi_n(q)$, $a-q^n$, $1-aq^n$, and $b-q^n$
are pairwise relatively prime, and all denominators are relatively
prime to their product. Combining these congruences after clearing
denominators proves \eqref{3.4}.

We may now set $b=1$. All specialized denominators remain relatively
prime to $\Phi_n(q)(a-q^n)(1-aq^n)$.
In particular, in the summands of $T_{n,s}$ we have
\[
k-s\le\frac{n-1}{2}+s<n,
\qquad
k\le\frac{n-1}{2}+2s\le n-2,
\]
so neither $(q^4;q^4)_{k-s}$ nor $(q^4;q^2)_k$
contains $\Phi_n(q)$.
The modulus in \eqref{3.4} becomes
$\Phi_n(q)(1-q^n)(a-q^n)(1-aq^n)$, which is divisible by
$\Phi_n(q)^2(a-q^n)(1-aq^n)$.
This proves \eqref{edu:mod all 3.3}.
\end{proof}

With the help of Propositions \ref{prop:mod-phia} and \ref{prop:mod all 3.3}, we are now in a position to prove Theorems \ref{thm:mod-phi3} and \ref{thm:analogue}.

{\noindent\it Proof of Theorem \ref{thm:mod-phi3}. }
Putting $ b\mapsto 1$ in Proposition \ref{prop:mod-phia}, we conclude that, modulo $  \Phi _{n}(q) (a-q^{n}) ( 1-aq^{n}) $,
\begin{align*}
T_{n,s}(a,1;q)
&\equiv \frac{ \left ( aq;q^{2} \right ) _{s} \left ( q/a;q^{2} \right ) _{-s} \left ( q^{-1};q^{4} \right ) _{s} \left ( q^{3};q^{4} \right ) _{\frac{n+1}{4}}  }{    \left ( aq^{5} ; q^{4}\right ) _{ s }\left ( q^{4} ; q^{4}\right ) _{\frac{n+1}{4} }}
\\
&\qquad \times \left [ \frac{1-aq^{n}}{1-a^{2}} \frac{  \left (  q^{5}/a ; q^{4}\right ) _{\frac{n+1}{4}} }{  \left ( q^{4}/a;q^{4} \right ) _{\frac{n+1}{4}}  }
			-\frac{a(a-q^{n})}{1-a^{2}}\frac{  \left (  aq^{5} ; q^{4}\right ) _{\frac{n+1}{4}} }{   \left ( aq^{4};q^{4} \right ) _{\frac{n+1}{4}}  } \right ]q^{s^{2}+3s}.
\end{align*}
We now specialize $ a\mapsto 1$. It is clear that $(1-q^{n})^{2}$ contains the factor $ \Phi _{n}(q)^{2}$.
By the L'H$\hat{o}$pital's rule, we have
		\begin{align}\label{3.7}
			&\lim_{a \to 1} \left [ \frac{1-aq^{n}}{1-a^{2}} \frac{  \left (  q^{5}/a ; q^{4}\right ) _{\frac{n+1}{4}} }{  \left ( q^{4}/a;q^{4} \right ) _{\frac{n+1}{4}}  }
			-\frac{a(a-q^{n})}{1-a^{2}}\frac{  \left (  aq^{5} ; q^{4}\right ) _{\frac{n+1}{4}} }{   \left ( aq^{4};q^{4} \right ) _{\frac{n+1}{4}} }  \right ]
\notag\\
&= \frac{  \left (  q^{5} ; q^{4}\right ) _{\frac{n+1}{4}} }{  \left ( q^{4};q^{4} \right ) _{\frac{n+1}{4}}  }\left[1-(1-q^n)\sum_{j=0}^{\frac{n-3}{4}}
\left(\frac{q^{4j+5}}{1-q^{4j+5}}-\frac{q^{4j+4}}{1-q^{4j+4}}\right)\right].	\end{align}
Therefore, modulo $ \Phi _{n}(q) ^{3}$,
	\begin{equation*}
		\begin{aligned}
			T_{n,s}(1,1;q)
&\equiv  \frac{ \left ( q;q^{2} \right ) _{s} \left ( q;q^{2} \right ) _{-s} \left ( q^{-1};q^{4} \right ) _{s} \left ( q^{3};q^{4} \right ) _{\frac{n+1}{4}} \left (  q^{5} ; q^{4}\right ) _{\frac{n+1}{4}} }{    \left ( q^{5} ; q^{4}\right ) _{ s } \left ( q^{4} ; q^{4}\right ) _{\frac{n+1}{4} } \left ( q^{4};q^{4} \right ) _{\frac{n+1}{4}}}q^{s^{2}+3s}
\\
&\qquad\times\left[1-(1-q^n)\sum_{j=0}^{\frac{n-3}{4}}
\left(\frac{q^{4j+5}}{1-q^{4j+5}}-\frac{q^{4j+4}}{1-q^{4j+4}}\right)\right]
.
		\end{aligned}
	\end{equation*}
Then, from the definition of $T_{n,s}(a,b;q)$ given in \eqref{def:T}, to prove Theorem \ref{thm:mod-phi3}, it suffices to 	
justify that
\begin{align}\label{cong:thm1-1}
&\frac{ \left ( q^{3};q^{4} \right ) _{\frac{n+1}{4}} \left (  q^{5} ; q^{4}\right ) _{\frac{n+1}{4}} }{    \left ( q^{4} ; q^{4}\right ) _{\frac{n+1}{4} } \left ( q^{4};q^{4} \right ) _{\frac{n+1}{4}}}
\left[1-(1-q^n)\sum_{j=0}^{\frac{n-3}{4}}
\left(\frac{q^{4j+5}}{1-q^{4j+5}}-\frac{q^{4j+4}}{1-q^{4j+4}}\right)\right]
\notag\\
&\qquad\equiv  -[n]q^{\frac{n+5}{4}  }\frac{\left(q^{-2} ; q^{4}\right)_{(n+1) / 4}}{\left(q^{4} ; q^{4}\right)_{(n+1) / 4}} \pmod{\Phi _{n}(q) ^{3}}.
\end{align}

Note that
\[
q^{-n}\equiv 1+(1-q^n)\pmod{\Phi_n(q)^2},
\]
since
\[
q^{-n}-1-(1-q^n)=\frac{(1-q^n)^2}{q^n}.
\]
Consequently, for $0\le j\le (n-3)/4$,
\[
\frac{1-q^{4j+5-n}}{1-q^{4j+4-n}}
\equiv
\frac{1-q^{4j+5}}{1-q^{4j+4}}
\left[
1-(1-q^n)
\left(
\frac{q^{4j+5}}{1-q^{4j+5}}
-\frac{q^{4j+4}}{1-q^{4j+4}}
\right)
\right]
\pmod{\Phi_n(q)^2}.
\]
All denominators occurring here are relatively prime to
$\Phi_n(q)$. Multiplying these congruences gives
\begin{align*}
\frac{(q^{5-n};q^4)_{\frac{n+1}{4}}}
     {(q^{4-n};q^4)_{\frac{n+1}{4}}}
\equiv
\frac{(q^5;q^4)_{\frac{n+1}{4}}}
     {(q^4;q^4)_{\frac{n+1}{4}}}
\left[
1-(1-q^n)
\sum_{j=0}^{\frac{n-3}{4}}
\left(
\frac{q^{4j+5}}{1-q^{4j+5}}
-\frac{q^{4j+4}}{1-q^{4j+4}}
\right)
\right]
\pmod{\Phi_n(q)^2}.
\end{align*}
Since $(q^3;q^4)_{\frac{n+1}{4}}$ contains the factor
$1-q^n$, it follows that
\begin{align}\label{cong:thm1-2}
&\frac{(q^3;q^4)_{\frac{n+1}{4}}
       (q^5;q^4)_{\frac{n+1}{4}}}
      {(q^4;q^4)_{\frac{n+1}{4}}^2}
\left[
1-(1-q^n)
\sum_{j=0}^{\frac{n-3}{4}}
\left(
\frac{q^{4j+5}}{1-q^{4j+5}}
-\frac{q^{4j+4}}{1-q^{4j+4}}
\right)
\right]
\notag\\
&\qquad\equiv
\frac{(q^3;q^4)_{\frac{n+1}{4}}
       (q^{5-n};q^4)_{\frac{n+1}{4}}}
      {(q^4;q^4)_{\frac{n+1}{4}}
       (q^{4-n};q^4)_{\frac{n+1}{4}}}
\pmod{\Phi_n(q)^3}.
\end{align}
Reversing the order of the factors, we obtain
\[
(q^{5-n};q^4)_{\frac{n+1}{4}}
=
(-1)^{\frac{n+1}{4}}
q^{-\frac{(n+1)(n-7)}8}
(q^{-2};q^4)_{\frac{n+1}{4}}
\]
and
\[
\frac{(q^3;q^4)_{\frac{n+1}{4}}}
     {(q^{4-n};q^4)_{\frac{n+1}{4}}}
=
(-1)^{\frac{n-3}{4}}
q^{\frac{(n-3)(n-1)}8}[n].
\]
Therefore,
\begin{equation}\label{eq:product-exact}
\frac{(q^3;q^4)_{\frac{n+1}{4}}
       (q^{5-n};q^4)_{\frac{n+1}{4}}}
      {(q^4;q^4)_{\frac{n+1}{4}}
       (q^{4-n};q^4)_{\frac{n+1}{4}}}
=
-[n]q^{\frac{n+5}{4}}
\frac{(q^{-2};q^4)_{\frac{n+1}{4}}}
     {(q^4;q^4)_{\frac{n+1}{4}}}.
\end{equation}
Substituting this into \eqref{cong:thm1-2}, we obtain \eqref{cong:thm1-1} as desired.
\qed
	
{\noindent\it Proof of Theorem \ref{thm:analogue}.}
We specialize $a=1$ in Proposition \ref{prop:mod all 3.3}.
This specialization is legitimate: the right-hand side has no pole
at $a=1$, and all specialized denominators are relatively prime to
$\Phi_n(q)$. For the left-hand side $T_{n,s}(a,b;q)$ given in \eqref{def:T}, in addition to the bounds used
in the proposition, we have
\[
k+s\le\frac{n-1}{2}+3s\le n-1.
\]
The denominators arising from negative indices are also relatively
prime to $\Phi_n(q)$, since their exponents have absolute value at
most $2s-1<n$.
The modulus becomes $\Phi_n(q)^2(1-q^n)^2$ and thus contains
$\Phi_n(q)^4$. By \eqref{eq:product-exact}, we obtain, modulo
$\Phi_n(q)^4$,
\begin{align}\label{eq:thm12-three-products}
T_{n,s}(1,1;q)
&\equiv
\frac{(q;q^2)_s(q;q^2)_{-s}(q^{-1};q^4)_s}
     {(q^5;q^4)_s}q^{s^2+3s}
\notag\\
&\quad\times\left\{
-[n]q^{\frac{n+5}{4}}
\frac{(q^{-2};q^4)_{\frac{n+1}{4}}}
     {(q^4;q^4)_{\frac{n+1}{4}}}
+\frac{(q^3,q^5;q^4)_{\frac{n+1}{4}}}
      {(q^4;q^4)_{\frac{n+1}{4}}^2}
-\frac{(q^3,q^5;q^4)_{\frac{n+1}{4}}}
      {(q^{4-n},q^{4+n};q^4)_{\frac{n+1}{4}}}
\right\}.
\end{align}
We now simplify the last two fractions, using a product expansion
as in the proof of Theorem \ref{thm:mod-phi3}.
For $1\le j\le(n+1)/4$, the exact identity
\[
(1-q^{4j-n})(1-q^{4j+n})
=(1-q^{4j})^2-q^{4j-n}(1-q^n)^2
\]
gives
\begin{align*}
\frac{(q^4;q^4)_{\frac{n+1}{4}}^2}
     {(q^{4-n},q^{4+n};q^4)_{\frac{n+1}{4}}}
&=\prod_{j=1}^{\frac{n+1}{4}}
\left(1-\frac{q^{4j-n}(1-q^n)^2}{(1-q^{4j})^2}\right)^{-1}
\\
&\equiv
1+(1-q^n)^2\sum_{j=1}^{\frac{n+1}{4}}
\frac{q^{4j}}{(1-q^{4j})^2}
\pmod{\Phi_n(q)^3}.
\end{align*}
Here terms containing $(1-q^n)^4$ have been discarded, and
$q^{-n}$ may be replaced by $1$ in the correction term because
it is already multiplied by $(1-q^n)^2$.
Since $(q^3;q^4)_{\frac{n+1}{4}}$ contains $1-q^n$, it follows that modulo $\Phi_n(q)^4$,
\begin{align}\label{eq:thm12-correction}
\frac{(q^3,q^5;q^4)_{\frac{n+1}{4}}}
      {(q^4;q^4)_{\frac{n+1}{4}}^2}
-\frac{(q^3,q^5;q^4)_{\frac{n+1}{4}}}
      {(q^{4-n},q^{4+n};q^4)_{\frac{n+1}{4}}}
\equiv
-[n]^2\frac{(q^3,q^5;q^4)_{\frac{n+1}{4}}}
                {(q^4;q^4)_{\frac{n+1}{4}}^2}
\sum_{j=1}^{\frac{n+1}{4}}\frac{q^{4j}}{[4j]^2}.
\end{align}
Moreover, reducing the shifted factors modulo $\Phi_n(q)$ and
again using the factor $1-q^n$ in $(q^3;q^4)_{\frac{n+1}{4}}$, we have
\begin{align*}
\frac{(q^3,q^5;q^4)_{\frac{n+1}{4}}}
     {(q^4;q^4)_{\frac{n+1}{4}}^2}
&\equiv
\frac{(q^3,q^{5-n};q^4)_{\frac{n+1}{4}}}
     {(q^4,q^{4-n};q^4)_{\frac{n+1}{4}}}
\\
&=-[n]q^{\frac{n+5}{4}}
\frac{(q^{-2};q^4)_{\frac{n+1}{4}}}
     {(q^4;q^4)_{\frac{n+1}{4}}}
\pmod{\Phi_n(q)^2}.
\end{align*}
Because the correction term in \eqref{eq:thm12-correction} contains
$[n]^2$, this replacement is valid modulo $\Phi_n(q)^4$.
Substituting into \eqref{eq:thm12-three-products}, we conclude that modulo $\Phi_n(q)^4$,
\begin{align*}
T_{n,s}(1,1;q)
\equiv -[n]q^{\frac{n+5}{4}}
\frac{(q^{-2};q^4)_{\frac{n+1}{4}}}
     {(q^4;q^4)_{\frac{n+1}{4}}}
\frac{(q;q^2)_s(q;q^2)_{-s}(q^{-1};q^4)_s}{(q^5;q^4)_s}
\left[1-[n]^2\sum_{j=1}^{\frac{n+1}{4}}
\frac{q^{4j}}{[4j]^2}\right]q^{s^2+3s},
\end{align*}
Combining the definition of $T_{n,s}(a,b;q)$ given in \eqref{def:T}, we finish the proof of Theorem \ref{thm:analogue}.
\qed	
	
\section{Proofs of Theorems \ref{thm.1.4} and \ref{thm.1.6}}\label{Sect.4}
To prove Theorems \ref{thm.1.4} and \ref{thm.1.6}, we first provide the following $q$-congruence.
\begin{thm}\label{thm.1.3}
Let $n> 1$ be an integer with $n \equiv 3\ (\bmod \ 4)$, and $s$ be an integer with $0\le s \le \frac{n-1}{6}$. Then, modulo $[n]$,
		\begin{equation*}
			\begin{aligned}
				\sum_{k=s}^{\frac{n-1}{2}+2s }& \left [ 6k+1 \right ] \frac{ \left ( q;q^{2} \right )_{k}  \left ( q;q^{2} \right )_{k-2s}  \left ( q;q^{2} \right )_{k+2s}  \left ( q^{-1};q^{4} \right )_{k}  }{ \left ( q^{4};q^{4} \right )_{k+s}  \left ( q^{4};q^{4} \right )_{k}  \left ( q^{4};q^{4} \right )_{k-s}  \left ( q^{4};q^{2} \right )_{k}  }q^{k^{2}+3k} \equiv 0 .
			\end{aligned}
		\end{equation*}
	\end{thm}
\begin{proof}
Replacing $k$ by $k+s$, it is equivalent to proving that
\begin{equation}\label{4.1}
\sum_{k=0}^{\frac{n-1}{2}+s }\omega _{q} \left ( k ;s\right ) \equiv 0 \pmod{[n]},
\end{equation}
where	
\begin{equation*}
\omega _{q} \left ( k;s \right ):=\left [ 6k+6s+1 \right ] \frac{ \left ( q;q^{2} \right )_{k+s}  \left ( q;q^{2} \right )_{k-s}  \left ( q;q^{2} \right )_{k+3s}  \left ( q^{-1};q^{4} \right )_{k+s}  }{ \left ( q^{4};q^{4} \right )_{k+2s}  \left ( q^{4};q^{4} \right )_{k+s}  \left ( q^{4};q^{4} \right )_{k}  \left ( q^{4};q^{2} \right )_{k+s}  }q^{(k+s)^{2}+3(k+s)}.
\end{equation*}
%
Let $\zeta $ be a primitive root of unity of degree $d$ with $d|n$ and $d>1$. It suffices to verify that for every divisor $d>1$ of $n$,
\begin{equation}\label{eq:all-d}
\sum_{k=0}^{\frac{n-1}{2}+s }\omega_q(k;s)
\equiv0\pmod{\Phi_d(q)}.
\end{equation}

Write
\[
s=td\pm r,\qquad t\in\mathbb Z_{\ge0},\qquad
0\le r\le\frac{d-1}{2}.
\]
We first identify conditions under which $\omega_q(k;s)$
tends to zero as $q\to\zeta$.
Set
\[
\rho_d:=\begin{cases}
(d+5)/4,&d\equiv3\pmod4,\\
(3d+5)/4,&d\equiv1\pmod4.
\end{cases}
\]
Thus $\rho_d$ is the smallest length for which $(q^{-1};q^4)_{\rho_d}$
contains a factor divisible by $\Phi_d(q)$.
Let $\nu_d$ count the factors $1-q^{jd}$ in the numerator minus those
in the denominator, where $j$ is a nonzero integer and repeated factors
are counted separately. Each such factor contains exactly one factor
$\Phi_d(q)$. For every integer $v$, including negative $v$, we have
\[
\nu_d\bigl((q;q^2)_v\bigr)
=\left\lfloor\frac{v+(d-1)/2}{d}\right\rfloor.
\]
Counting the corresponding factors in
\[
F_q(k;s):=\frac{\omega_q(k;s)}{[6k+6s+1]},
\]
and writing $u$ for
the least nonnegative residue of $k+s$ modulo $d$, gives
\begin{align*}
\nu_d(F_q(k;s))=&
\left\lfloor\frac{u+(d-1)/2}{d}\right\rfloor
+\left\lfloor\frac{u-2r+(d-1)/2}{d}\right\rfloor
+\left\lfloor\frac{u+2r+(d-1)/2}{d}\right\rfloor\\
&+\left\lfloor\frac{u+d-\rho_d}{d}\right\rfloor
-\left\lfloor\frac{u+r}{d}\right\rfloor
-\left\lfloor\frac ud\right\rfloor
-\left\lfloor\frac{u-r}{d}\right\rfloor
-\left\lfloor\frac{u+1}{d}\right\rfloor.
\end{align*}
The formula holds for both signs: replacing $r$ by $-r$ interchanges
the terms involving $u-2r$ and $u+2r$, and those involving $u-r$ and $u+r$.
Evaluating the floors shows that the expression is nonnegative.
It is zero precisely when
\[
r\le u\le\frac{d-1}{2}-2r,\qquad u<\rho_d,
\]
apart from the exceptional case $(d,r,u)=(5,1,4)$.
In that case, $[6k+6s+1]$ is divisible by $\Phi_5(q)$,
since $6u+1=25$.

Thus all the limits
\[
\omega_\zeta(k;s):=\lim_{q\to\zeta}\omega_q(k;s)
\]
exist, and a nonzero limit is possible only if
\begin{equation}\label{eq:possible-support}
r\le u\le\frac{d-1}{2}-2r.
\end{equation}
Indeed, if $r>(d-1)/6$, this interval is empty.
If $r=(d-1)/6$, its only possible value is $u=r$,
and $[6k+6s+1]$ vanishes at $q=\zeta$
because $6r+1=d$.
Consequently, all the limits vanish when
$r\ge(d-1)/6$.
The desired congruence \eqref{eq:all-d} follows immediately in this case.

It remains to consider $0\le r<(d-1)/6$.
Applying Lemma \ref{lem 1} with $n=d$, $s=r$ and $a=b=1$, and then replacing $k$ by $k+r$, gives
\[\sum_{k=0}^{\frac{d-1}{2}+r}\omega_\zeta(k;r)=0.\]
By \eqref{eq:possible-support},
$\omega_\zeta(k;r)=0$ for
$k>(d-1)/2-3r$.
Therefore,
\begin{equation}\label{eq:basic-zero}
\sum_{k=0}^{d-1}\omega_\zeta(k;r)
=\sum_{k=0}^{\frac{d-1}{2}}\omega_\zeta(k;r)
=\sum_{k=0}^{\frac{d-1}{2}-3r}\omega_\zeta(k;r)=0.
\end{equation}
If $s=td+r$ with $0\le r<(d-1)/6$, $\omega_\zeta(\ell d;s)$ and $\omega_\zeta(0;r)$ are nonzero. Considering the following relation
\begin{equation*}
\frac{\omega _{\zeta } \left ( ld+k;s \right )}{\omega _{\zeta} \left ( ld;s \right )} = \lim_{q \to \zeta } \frac{\omega _{q} \left ( ld+k;s \right )}{\omega _{q} \left ( ld;s \right )} =\frac{\omega _{\zeta} \left ( k;r\right )}{\omega _{\zeta} \left ( 0;r \right )},
\end{equation*}
it follows from \eqref{eq:basic-zero} that
\[
\sum_{k=0}^{d-1}\omega_\zeta(\ell d+k;s)
=
\frac{\omega_\zeta(\ell d;s)}{\omega_\zeta(0;r)}
\sum_{k=0}^{d-1}\omega_\zeta(k;r)
=0.
\]
If $s=td-r$, both $\omega_\zeta(\ell d+2r;s)$ and
$\omega_\zeta(0;r)$ are nonzero. We obtain
\[
\frac{\omega_\zeta(\ell d+2r+k;s)}
     {\omega_\zeta(\ell d+2r;s)}
=
\lim_{q\to\zeta}
\frac{\omega_q(\ell d+2r+k;s)}
     {\omega_q(\ell d+2r;s)}
=
\frac{\omega_\zeta(k;r)}{\omega_\zeta(0;r)}.
\]
Consequently, using \eqref{eq:possible-support} and \eqref{eq:basic-zero} gives
\begin{align*}
\sum_{k=0}^{d-1}\omega_\zeta(\ell d+k;s)
&=\sum_{k=2r}^{(d-1)/2-r}\omega_\zeta(\ell d+k;s)=
\sum_{k=0}^{(d-1)/2-3r}
\omega_\zeta(\ell d+2r+k;s)\\
&=
\frac{\omega_\zeta(\ell d+2r;s)}
     {\omega_\zeta(0;r)}
\sum_{k=0}^{(d-1)/2-3r}\omega_\zeta(k;r)
=0.
\end{align*}

Furthermore, the truncated block with upper limit
$(d-1)/2\pm r$ contains all its possibly nonzero terms.
Hence
\[
\sum_{k=0}^{(d-1)/2\pm r}
\omega_\zeta(\ell d+k;s)=\sum_{k=0}^{d-1}\omega_\zeta(\ell d+k;s)=0,
\]
where the sign agrees with that in $s=td\pm r$.

Finally, since
\[
\frac{n-1}{2}+s
=
\left(\frac{n-d}{2d}+t\right)d
+\frac{d-1}{2}\pm r,
\]
we obtain
\begin{align*}
\sum_{k=0}^{(n-1)/2+s}\omega_\zeta(k;s)
&=
\sum_{\ell=0}^{(n-d)/(2d)+t-1}
\sum_{k=0}^{d-1}\omega_\zeta(\ell d+k;s)+
\sum_{k=0}^{(d-1)/2\pm r}
\omega_\zeta\left(\frac{n-d}{2}+td+k;s\right)
=0,
\end{align*}
which yields that \eqref{eq:all-d} is valid. This completes the proof.
\end{proof}	

According to Theorem \ref{thm.1.3}, we can finish the proof of Theorems \ref{thm.1.4} and \ref{thm.1.6}.

{\noindent\it Proof of Theorem \ref{thm.1.4}.}	
Combining Theorem \ref{thm:mod-phi3} and Theorem \ref{thm.1.3}, the proof follows from the fact that the least common multiple of $[n]$ and $\Phi _{n}(q)^{3} $ is $[n] \Phi _{n}(q)^{2}$.
\qed
	
{\noindent\it Proof of Theorem \ref{thm.1.6}.} Since the least common multiple of $[n]$ and $\Phi _{n}(q)^{4} $ is $[n] \Phi _{n}(q)^{3}$, we complete the proof by using Theorem \ref{thm:analogue} and Theorem \ref{thm.1.3}.
\qed

\vspace{0.5cm}
 \baselineskip 15pt
{\noindent\bf\large{\ Acknowledgements}} \vspace{7pt} \par
The authors would like to thank their supervisor, Helen W. J. Zhang, for continuous guidance and support.

\end{document}